\documentclass[11pt]{amsart}
\usepackage{amssymb,amsthm,amsmath}
\RequirePackage[dvipsnames,usenames]{xcolor}
\usepackage{hyperref}
\usepackage{mathtools,etoolbox}
\usepackage[all]{xy}
\usepackage{tikz}
\usepackage{enumitem}
\usepackage[english]{babel}
\usepackage{graphics}
\usepackage{color}

\hypersetup{
bookmarksdepth=3,
bookmarksopen,
bookmarksnumbered,
pdfstartview=FitH,
colorlinks,%backref,hyperindex,
linkcolor=Sepia,
anchorcolor=BurntOrange,
citecolor=MidnightBlue,
citecolor=OliveGreen,
filecolor=BlueViolet,
menucolor=Yellow,
urlcolor=OliveGreen
}

\newcommand{\factor}[2]{\left. \raise 2pt\hbox{\ensuremath{#1}} \right/
        \hskip -2pt\raise -2pt\hbox{\ensuremath{#2}}}

\makeatletter
\renewcommand\subsection{
  \renewcommand{\sfdefault}{pag}
  \@startsection{subsection}%
  {2}{0pt}{.8\baselineskip}{.4\baselineskip}{\raggedright
    \sffamily\itshape\small\bfseries
  }}
\renewcommand\section{
  \renewcommand{\sfdefault}{phv}
  \@startsection{section} %
  {1}{0pt}{\baselineskip}{.8\baselineskip}{\centering
    \sffamily
    \scshape
    \bfseries
}}
\makeatother

\usepackage[left=1.02in,top=1.0in,right=1.02in,bottom=1.0in]{geometry}
\usepackage{multirow}

\setlist[enumerate]{leftmargin=0.8cm}
\setlist[itemize]{leftmargin=0.8cm}
\setlist[description]{leftmargin=0.0cm}

\newtheorem{theorem}{Theorem}[section]

\newtheorem{proposition}[theorem]{Proposition}
\newtheorem{conjecture}[theorem]{Conjecture}
\newtheorem{lemma}[theorem]{Lemma}

\theoremstyle{definition}

\newtheorem{example}[theorem]{Example}
\theoremstyle{remark}

\usepackage{color}

\DeclareMathOperator{\Exc}{Exc}

\numberwithin{equation}{section}

\newcommand{\A}{{\mathbb A}}
\newcommand{\G}{{\mathbb G}}

\newcommand{\PP}{{\mathbb P}}

\newcommand{\C}{{\mathbb C}}

\newcommand{\kbar}{{\overline{k}}}

\newcommand{\calC}{{\mathcal C}}

\newcommand{\calO}{{\mathcal O}}

\usepackage[mathscr]{euscript}

\newcommand{\scrD}{{\mathscr D}}

\newcommand{\scrX}{{\mathscr X}}

\DeclareMathOperator{\Spec}{Spec}

\numberwithin{equation}{section}
\numberwithin{table}{section}

\title{Some examples of exceptional loci in Vojta Conjecture}
\author{Amos Turchet}
\address{Dipartimento di Matematica e Fisica, Università degli Studi Roma 3, L.go S. L. Murialdo 1, 00146 Roma, Italy}
\email{amos.turchet@uniroma3.it}

\subjclass[2010]{14G40, 14G05.}
\keywords{Vojta conjecture, exceptional locus, integral points, hyperbolicity.}

\begin{document}
\begin{abstract}
  In this short note we discuss the exceptional locus for the Lang-Vojta's conjecture in the case of the complement of two completely reducible hyperplane sections in a cubic surface. Using elementary methods, we show that generically the exceptional set is the union of the remaining 21 lines in the surface. We also describe examples in which the exceptional set is strictly larger.
\end{abstract}
\maketitle
%\tableofcontents
\section{Introduction}\label{sec:intro}

One of the main questions in Diophantine Geometry is to determine whether, given an algebraic variety $X$ defined over a number field $k$, the closure of the set of $k$-rational points $\overline{X(k)}$ is a proper subset of $X$. A refined and strictly stronger question asks if there exists a proper closed subset $Z$ of $X$ (called the \emph{exceptional set}) such that, for every finite extension $k' \supset k$, the complement $X \setminus Z$ contains only finitely many $k'$-rational points. Similar questions can be formulated for quasi-projective varieties and $S$-integral points.

Conjectures of Lang, Vojta and Campana, as in \cite[Conjecture 3.7]{Lang91}, \cite[Conjecture 3.4.3]{Vojta87} and \cite[Conjectures 9.2 and 9.20]{Ca04} (see \cite{BG,SH,AT_book,DT_invit} for introductions and discussions), predict that, if $X$ is of general (resp. of log general) type there exists an exceptional set and hence the rational (resp. integral) points are not Zariski dense. In fact, the stronger versions of such conjectures predict that the existence of an exceptional set should be equivalent to the fact that $X$ is of general (resp. log general) type.

Despite several deep results of various authors (for example \cite{Falt_Ab,VojtaSA,McQ98,VojtaSA2,CZAnnals,Levin,CoZa10,Aut2,RV19} and many others), all the above conjectures are still widely open in general, at least when the dimension of $X$ is greater than 1. 

The main goal of this paper is to discuss the exceptional set in the special case where $X$ is the complement of a specific divisor in a smooth cubic surface. In this setting, \cite[Theorem 1]{CoZa10} shows that $S$-integral points on $X$, for a finite set of places $S$, are not Zariski dense. Here we compute the exceptional set and discuss other related questions.

The main result is the following: 

\begin{theorem}
  \label{thm:cubic}
  Let $X$ be a generic smooth cubic surface defined over $k$ and let $D$ be the divisor consisting of two completely reducible hyperplane sections of $X$. Then the exceptional locus $\Exc(X,D)$ consists of the 21 lines in $X_{\kbar}$ not contained in $D$.
\end{theorem}

%The statement of Theorem \ref{thm:cubic} appears in \cite{CoZa10} without proof.\Amos{rephrase}

We show that the genericity hypothesis in Theorem \ref{thm:cubic} cannot be removed by presenting some explicit examples where the exceptional set strictly contains the 21 lines outside of $D$.

\subsection*{Acknowledgements.} This project was part of the author Master Thesis under the supervision of Pietro Corvaja. We thank him for his support and several discussions. We thank Barbara Bolognese and Laura Capuano for comments and discussions. The author is partially supported by the project PRIN2017: Advances in Moduli Theory and Birational Classification and he is a member of the GNSAGA INdAM group.

\section{The cubic surface}
\label{sec:cubic}

We collect in this section standard facts on cubic surfaces.%that will be used throughout the article.

Let $X$ be a smooth cubic surface defined over an algebraically closed field of characteristic 0. Then, $X$ is isomorphic to the blow-up of $\PP^2$ at six points $P_1,\dots,P_6$ in general position \cite[V, Corollary 4.7]{Hartshorne}. It is a fundamental result, originally proved by Cayley and Salmon, that every such surface contains exactly 27 lines. More precisely, the following holds:

\begin{theorem}[{\cite[V, Theorem 4.9]{Hartshorne} and \cite[Lemma IV,5]{Beauville}}]
  \label{thm:27}
  Let $X,P_1,\dots,P_6$ as above. Then, there are exactly 27 lines in $X$, each one of self intersection -1, and they are:
  \begin{itemize}
    \item the 6 exceptional lines $E_i$ of the blow up lying over $P_i$, for $i=1,\dots,6$;
    \item the 15 strict transforms of the lines $\ell_{i,j}$ passing through $P_i,P_j$ for $i,j \in \{1,\dots,6\}$, $i < j$;
    \item the 6 strict transforms of the conics $C_i$ passing through $\{ P_1, \dots, P_6\} \setminus \{ P_i\}$ for $i = 1,\dots,6$.
  \end{itemize}
  Moreover, given any line $L$ in $X$, there are exactly 10 other lines (distinct from $L$) that intersects $L$, and fall into 5 disjoint pairs of concurrent lines.  
\end{theorem}

Given any smooth cubic surface $X$, any choice of 6 skew lines in $X$ allows to construct a morphism $X \to \PP^2$ that exhibits $X$ as the blow up of six points in general position and such that the chosen 6 lines are the exceptional lines of the blow up \cite[V, Proposition 4.10]{Hartshorne}. 

Using this construction we will make the following choice for the hyperplane sections.\medskip

\paragraph{\textbf{Setting}}
Let $X$ be a smooth cubic surface  and let $D = H_1 + H_2$ be the divisor corresponding to two completely reducible hyperplane sections, both defined over a number field $k$. Let $\pi: X \to \PP^2$ be the morphism that exhibits $X$ as the blow up of 6 points in general position. Without loss of generality, using the notation of Theorem \ref{thm:27}, we can assume that $H_1,H_2$ correspond to the strict transforms of $\ell_{1,2} + \ell_{3,4} + \ell_{5,6}$, and $\ell_{2,3} + \ell_{4,5} + \ell_{1,6}$ respectively, i.e.
\begin{equation}
  \label{eq:H}
  H_1 = \pi^{-1}_* \left( \ell_{1,2} + \ell_{3,4} + \ell_{5,6} \right) \qquad \qquad H_2 = \pi^{-1}_* \left( \ell_{2,3} + \ell_{4,5} + \ell_{1,6} \right).
\end{equation}

From now, on we will consider $X,D = H_1 + H_2$ fixed as in the above setting.

\section{The exceptional set}\label{sec:exc}
Let $(X,D)$ be a pair, where $X$ is a normal projective variety and $D$ is a normal crossing divisor. The pair is said of \emph{log general type} if for a(ny) log resolution $(Y,D_Y)$ the divisor $K_Y + D_Y$ is big. From the arithmetic point of view, a pair $(X,D)$ defined over a number field $k$ encodes the quasi-projective variety $X \setminus D$ and is the natural ambient to study questions about $(S,D)$-integral points.
In particular, given a model $(\scrX,\scrD)$ over $\Spec \calO_k$ of the pair $(X,D)$, the set of $(S,D)$\emph{-integral points} is the set of sections $\sigma: \Spec \calO_k \to \scrX$ such that the support of $\sigma^* \scrD$ is contained in $S$. If $(X,D)$ is given inside a projective space, the $(S,D)$-integral points correpond to $k$-rational points in $X \setminus D$ whose reduction modulo every prime outside $S$ is not contained in the reduction of $D$. In particular, if $D$ is empty one recovers the notion of $k$-rational points. For more details and standard properties we refer to \cite[Section 3]{AT_book}.

One of the fundamental conjectures that drives the study of rational and integral points on varieties is the following due to Lang and Vojta.

\begin{conjecture}[Strong Lang-Vojta]
  \label{conj:LV}
Let $(X,D)$ be a pair defined over a number field $k$. Then $(X,D)$ is of log general type if and only if there exists a proper closed subset $Z = \Exc(X,D)$, called the \emph{exceptional set}, such that $X \setminus Z$ has only finitely many $(S',D)$-integral points, for every finite extension $k' \supset k$, where $S'$ is the set of places above $S$.
\end{conjecture}

The conjecture has a series of far reaching implications, and it is known only in some particular cases: we refer to \cite{AT_book,advt,abt,J_book} for examples and discussions.
The exceptional set $\Exc(X,D)$ in this article is sometimes called the \emph{Diophantine} exceptional set to distinguish it from the \emph{holomorphic} exceptional set (i.e. the union of all images of entire curves in the analytification $X_\C$). These sets were introduced by Serge Lang in \cite{Lan86} and he conjectured that the two sets coincide. For a detailed discussion about exceptional sets we refer to \cite[Section 3]{Levin} and \cite[Section 3]{RTW} for some generalization.

In our setting, where $X$ is a smooth cubic surface and $D = H_1 + H_2$, the divisor $K_X + D$ is ample so in particular the pair $(X,D)$ is of log general type and Conjecture \ref{conj:LV} predicts that $\Exc(X,D)$ is a proper closed subset of $X$.  In \cite[Theorem 1]{CoZa10} Corvaja and Zannier showed that the $\calO_{S,D}$ integral points on $X$ are not Zariski dense; in fact, using for example \cite{Levin}, or directly \cite[Theorem 1.6]{RTW2}, the exceptional locus is a proper subvariety of $X$.

Note that by definition, since $X$ is a surface, $\Exc(X,D)$ contains any quasi-projective curve $\calC \setminus D_{\calC}$ contained in $X \setminus D$, such that $\calC(\calO_{S,D_{\calC}})$ is potentially infinite, i.e. infinite up to possibly a finite extension of the base field. By Siegel's Theorem, the normalization of every such curve is either a (projective) elliptic curve, or the complement of at most two points in a rational curve. 
Since $H_1$ and $H_2$ are hyperplane sections, every projective curve $\calC$ intersects each of them in at least one point. The above discussion implies the following characterization of the exceptional set.

\begin{lemma}
  \label{lem:exc}
  Let $(X,D)$ be as in our setting. Then, the exceptional set $\Exc(X,D)$ satisfies:
  \[
   \Exc(X,D) \subset \{ \calC \mid \calC^\nu \cong \PP^1\subset X \text{ and }\# \nu^* D \leq 2 \} = \{ f(\PP^1) \mid f\colon \PP^1 \to X \text{ and } \# f^*D \leq 2 \},
    \]
    where $\nu: \calC^\nu \to \calC$ is the normalization.
\end{lemma}

Note that the cardinalities in the above formula are the cardinalities of the \emph{support} of the divisors $\nu^*D$ and $f^*D$, disregarding multiplicities. In other words, we are interested in maps $f: \PP^1 \to X$ such that the truncated counting function $N^{[1]}_f(D)$ (see \cite[Section 22]{Vojta11} for details)  satisfies $N^{[1]}_f(D) \leq 2$ . In particular, the intersection between $f(\PP^1)$ and $D$ can have high multiplicity.

\section{Main result}
In this section we are going to prove Theorem \ref{thm:cubic}, namely that for a \emph{generic} cubic surface $X$, the excpetional set $\Exc(X,D)$ is the union of all the lines not contained in $D = H_1 + H_2$.

We will begin by proving that the image of every morphism $f: \PP^1 \to X$ such that the support of $f^*D$ consists of two points has degree bounded by 2. This is related to the fact that fact that $(X,D)$ is algebraically hyperbolic \cite[Theorem 8.2]{RTW} as predicted by the analogue of Conjecture \ref{conj:LV} (we refer to \cite{CZConic,CZGm,Tur,CaTur,JK,RTW2,GSW} for examples and definitions).

\begin{proposition}\label{prop:deg_bound}
Let $(X,D)$ be as in our setting and assume that $X$ is generic (in the moduli space of cubic surfaces). Let $f: \PP^1 \to X$ be a morphism such that $f^*D$ is supported in at most 2 points. Then, either $f(\PP^1)$ is contracted by $\pi$ or
\[
  \deg \pi(f(\PP^1)) \leq 2.
\]
\end{proposition}

\begin{proof}
  Since $D$ is an ample divisor in $X$, the intersection with the image $f(\PP^1)$ is not empty. Denote by $Q_1, Q_2$ the points in $\pi(f(\PP^1) \cap D)$ (with possibly $Q_1 = Q_2$), where we recall that $\pi: X \to \PP^2$ is the blow up of $\PP^2$ at six points. We denote by $P_1,\dots,P_6$ the six points in general position that are the center of $\pi$. In our setting $\pi(D)$ looks as in Figure \ref{fig:6lines}.
\begin{figure}[h]
  \begin{centering}
    \includegraphics[scale=.5]{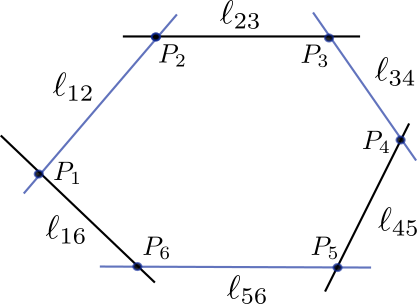}
    \caption{A graphical visualization of $\pi(D)$.}
    \label{fig:6lines}
  \end{centering}
\end{figure}

    Denote by $\calC = \pi(f(\PP^1))$. If $\calC$ is a point, then $f(\PP^1)$ is one of the exceptional lines of $\pi$ which proves the assertion. We can therefore assume that $\calC$ is a curve in $\PP^2$.

    The genericity hypothesis implies that we can assume, up to reordering, that there exists an index $j$, with $2 \leq j \leq 6$, such that both $Q_1$ and $Q_2$ do not lie neither on $\ell_{12}$ nor on $\ell_{j,j+ 1}$ (where if $j = 6$ we consider $j+1 = 1$). If we denote by $m_i$ the multiplicity of $\calC$ at $P_i$ then the degree of $\calC$ is equal to $m_1 + m_2 = m_{j} + m_{j+ 1}$, since the intersection has to be transverse.

%    We claim that there exists $j,k$ different from $1$ and $\ihat$ such that the intersection cycles are
%    \begin{align*}
%      &\calC \cdot \ell_{j, j+1} = m_j P_j + m_{j+1} P_{j+1} \\
%      &\calC \cdot \ell_{k, k+1} = m_k P_k + m_{k+1} P_{k+1}.
%    \end{align*}
%  This is true if there exists $j,k$ as above with $Q_1,Q_2 \notin \ell_{j,j+1},\ell_{k,k+1}$. If this is not the case then $Q_1$ is one of the $P_i$'s; let $i$ be such that $Q_1 = P_i$. Then $\calC$ is tangent to at most one between $\ell_{i-1,i}$ and $\ell_{i,i+1}$, since otherwise $Q_1 = Q_2$ and one can find two indices $j,k$ as wanted. Then we can find $j,k$, different from $1$ and $\ihat$ such that $\calC$ is not tangent to $l_{j,j+1}$ and $l_{k,k+1}$, and such that $\calC$ intersects both lines only in the points $P_j,P_{j+1},P_{k},P_{k+1}$. Which proves the claim.
%
%  It follows that without loss of generality the intersection of $\calC$ with 4 of the 6 lines is as follows, with $j < k$:
%  \begin{align}\label{eq:int_4}
%      &\calC \cdot \ell_{12} = m_1 P_1 + m_{2} P_{2} \nonumber \\
%      &\calC \cdot \ell_{23} = m_2 P_2 + m_{3} P_{3} \nonumber \\
%      &\calC \cdot \ell_{j, j+1} = m_j P_j + m_{j+1} P_{j+1} \\
%      &\calC \cdot \ell_{k, k+1} = m_k P_k + m_{k+1} P_{k+1}. \nonumber
%  \end{align}

%The sum of the coefficients in every line of \eqref{eq:int_4} is equal to the degree of $\calC$. It follows that $m_1 = m_3$ and there are three possible cases: either $m_2 = m_4$ (when $j=3$), $m_2 = m_6$ (when $k=6$) or $m_4 = m_6$ (when $j=4$ and $k=5$). 
Denote by $d = \deg \calC$, so that $d= m_1 + m_2 = m_{j} + m_{j+ 1}$. Assume that $j\neq 2,6$. Consider now the line $\ell_{1,j}$ that passes through $P_1,P_{j}$. The intersection $\calC \cdot \ell_{1,j}$ contains the cycle $m_1 P_1 + m_{j} P_{j}$ so that $d \geq m_1 + m_{j}$. 
This implies that $m_2 \geq m_{j}$. Using similarly the lines $\ell_{1,j+ 1}$, $\ell_{2,j}$ and $\ell_{2,j+ 1}$ one gets $m_1 = m_2 = m_{j} = m_{j+ 1}$.  

By \cite[Theorem 4.4]{Walker} one has that
\begin{equation}\label{eq:deg}
  (d-1)(d-2) \geq \sum_{i=1}^6 m_i(m_i-1).
\end{equation}
 
Expanding the formula with $d = m_1 + m_2$ we get
\begin{equation*}\label{eq:deg2}
    -(m_1- m_2)^2 - (m_1 + m_2) + 2 \geq 0,
\end{equation*}
which implies that $d = m_1 + m_2 \leq 2$ as wanted. 

It remains to deal with the case $j = 6$ (which is equivalent to the case $j=2$). In this case as above we get $m_3 = m_5$. Applying the same argument one gets $d \leq 4$, or more precisely $d \leq 2$ a part from the following two cases: either $m_1 = m_2 = m_6 = 2$ and all other multiplicities are equal to 1, or $m_1 = 2$ and all the other multiplicities are equal to 1. The first case corresponds to a quartic curve with three nodes at $P_1,P_2,P_6$, passing through the six points and through $\ell_{23} \cap \ell_{45}$ and $\ell_{34} \cap \ell_{56}$, and tangent at these points to two of the four lines. For a generic choice of the base points there will be no rational nodal quartic curve with these properties (see Example \ref{ex:quartic}). The second case corresponds to a cubic passing through the same eight points with a node in $P_1$ which again is excluded by genericity of the points (see Examples \ref{ex:cubic} and \ref{ex:cubic2}). This shows again that $d \leq 2$ as wanted.

%We will first deal with the case $m_2 = m_4$ (the case $m_2 = m_6$ being completely analogue).
%
%\begin{description}
%  \item[$\boxed{m_2 = m_4}$] Using the formula \eqref{eq:deg} with $d = m_1 + m_2$ we get
%  \begin{equation*}\label{eq:deg2}
%    -(m_1- m_2)^2 - (m_1 + m_2) + 2 \geq m_5^2 - m_5 + m_6^2 - m_6 \geq 0
%  \end{equation*}
%  which implies that $d = m_1 + m_2 \leq 2$ as wanted.
%  \item[$\boxed{m_4 = m_6}$] Consider now the line $\ell_{14}$ that passes through $P_1,P_4$. The intersection $\calC \cdot \ell_{1,4}$ contains the cycle $m_1 P_1 + m_4 P_4$ so that $d \geq m_1 + m_4$. This implies that $m_2 \geq m_4$. Using similarly the line $\ell_{2,5}$ one gets $m_1 \geq m_5$. But since $d = m_1 + m_2 = m_4 + m_5$ we get that $m_1 = m_5 $ and $m_2 = m_4$. Therefore we can apply the previous case and this concludes the proof.
%\end{description}

\end{proof}

We use the above proposition to show that $\calC$ is either a line passing through two of the points $P_i$ or a conic passing through 5 of the points.
\begin{proof}[Proof of Theorem \ref{thm:cubic}]
By Lemma \ref{lem:exc} the exceptional set $\Exc(X,D)$ is contained in the image of all morphisms $f: \PP^1 \to X$ such that $\# f^*D \leq 2$. Moreover, by Proposition \ref{prop:deg_bound}, the image $\calC = \pi(f(\PP^1))$ is either a point, in which case $f(\PP^1)$ is one of the six exceptional lines, or a plane curve of degree at most 2. Note that, by definition, the intersection between $\calC$ and $\pi(D)$ is supported in a set of the form $\{ P_1,\dots,P_6,Q_1,Q_2 \}$, where the points $Q_i$ might not be distinct from the points $P_j$. We distinguish three cases:

\begin{description}
  \item[$\boxed{Q_1,Q_2 \notin \{ P_1, \dots, P_6 \}}$] In this case the intersection of $\calC$ with $\pi(D)$ is transverse at $P_1,\dots,P_6$, and we can assume that $Q_1 \in \ell_{12}$. If $\calC$ is a line, then it has to intersect at least 4 of the lines $\ell_{i,i+1}$ in the points $P_3, P_4, P_5$ and $P_6$. The general position hypothesis implies that at most two of them lie on a line which shows that $\calC$ is a line through two of them. 

  Assume that $\calC$ is a conic. If $Q_2 \notin \ell_{23}$ or $Q_2 \notin \ell_{6,1}$ then it is easy to see that $P_1,P_2 \in \calC$ and thus the intersection with $\ell_{12}$ would contain three distinct points which is impossible. Therefore we can assume that $Q_2 \in \ell_{23}$ (the case $Q_2 \in \ell_{6,1}$ beign completely analogous). We claim that in this case $\calC$ is the unique conic passing through $P_1,P_3,P_4,P_5$ and $P_6$. This follows from the fact that $\calC$ has to pass thrugh $P_1$ and $P_3,\dots,P_6$ since the interesection cycle with each of the lines $\ell_{j,j+1}$, with $j=3,\dots,6$ has degree two and all the intersections are transverse at smooth points of $\calC$. 
  \item[$\boxed{Q_1 \in \{P_1,\dots,P_6 \}}$] We can assume $Q_1 = P_1$, and therefore $\calC$ is tangent in $P_1$ at one of the lines, which we can assume to be $\ell_{12}$. %, and intersects transversally $\pi(D)$ at all the points $P_2,\dots,P_6$. 
  In particular $\calC$ is not a line.
  
  This tangency condition implies that $P_2 \notin \calC$, hence that $Q_2 \in \ell_{23}$. If $\calC$ is tangent to $\ell_{23}$ in $Q_2$ then $P_3 \notin \calC$ which is impossible since $\calC$ would intersect $\ell_{34}$ transverslly at a single point. Therefore, $\calC$ is not tangent to $\ell_{23}$ in $Q_2$, which, with the same argument as in the previous case, implies that $P_1,P_3,P_4,P_5$ and $P_6$ are all points of $\calC$ and hence $\calC$ is the unique conic through these points.

  \item[$\boxed{Q_1,Q_2 \in \{P_1,\dots,P_6 \}}$] This in particular implies that $\calC$ does not intersect $\pi(D)$ transversally at $Q_1$ and $Q_2$ and therefore $\calC$ cannot be a line. 
  
  Then, $\calC$ is a conic and it is tangent to two of the lines $\ell_{j,j+1}$ at the points $Q_1$ and $Q_2$. This will force the curve to pass through at least 5 of the points as wanted. We note that the genericity hypothesis implies that there will be no such curves, since the points $P_i$ are generic.
\end{description}\end{proof}

\section{The non-generic case}
In the proofs of Proposition \ref{prop:deg_bound} and Theorem \ref{thm:cubic} we used the fact that the surface $X$ was generic to show that the triples of lines $\ell_{12},\ell_{3,4},\ell_{56}$ and $\ell_{23},\ell_{45}, \ell_{16}$ do not share a common point. More precisely, in the proofs we used that the points $P_1,\dots,P_6$ and the intersection of any two lines connecting two of them are in general position. We collect here some examples that show that, when this is not the case, the exceptional set has different properties. 

\begin{example}[Affine lines]\label{ex:lines}
 Theorem \ref{thm:cubic} shows that in the exceptional set the curves are all isomorphic to $\G_m$. For some special configuration of the base points one can show that some of the curves will be isomorphic to the affine line $\A^1$, so in particular they intersect $D$ transversally at only one point.
 
 Let $R = \ell_{23} \cap \ell_{56}$ and assume that $R \in \ell_{14}$, i.e. $R, P_1$ and $P_4$ are collinear. Then $\ell_{14}$ intersects each of the lines in $D$ transversally in exactly one point, and its strict transform intersects $D$ only in the (preimage of the) point $R$. In particular $\pi_*^{-1}(\ell_{14} \setminus \{ R \}) \subset \Exc(X,D)$ and $\ell_{14} \setminus \{ R \} \cong \A^1$.
\end{example}

\begin{example}[Extra cubics]\label{ex:cubic}
  Assume that there exist points $R_1 = \ell_{12} \cap \ell_{3,4} \cap \ell_{56}$ and $R_2 = \ell_{23} \cap \ell_{45} \cap \ell_{16}$. Then, any cubic passing throuhg $P_1,\dots,P_6$ transversally and through $R_1,R_2$ intersects $\pi(D)$ in only these points. In particular, its strict transform will intersect $D$ in only two points. Every nodal cubic through these points will then be contained in the exceptional set (and in general there will be 12 such singular cubics).
\end{example}

\begin{example}[Extra cubics 2]\label{ex:cubic2}
  Let $R_1 = \ell_{23} \cap \ell_{45}$ and let $R_2 = \ell_{34} \cap \ell_{56}$. For a generic choice of the points $P_1,\dots,P_6$ there will be no nodal cubic passing through $P_2,\dots,P_6,R_1,R_2$ and with a node in $P_1$ tangent to neither $\ell_{12}$ nor $\ell_{16}$. On the other hand, there are special choices of the points for which this will be possible. For example one could fix all the points but $P_2$, pick a nodal cubic with node in $P_1$ and passing through $P_3,\dots,P_6,R_1,R_2$ and then choose a point $P_2$ in the cubic such that $P_1,\dots,P_6$ are in general position. In this case the strict transform of the nodal cubic will be a rational curve in $X$ that intersects $D$ in only two points.
\end{example}

\begin{example}[Extra quartic]\label{ex:quartic}
 We have already seen in the proof of Proposition \ref{prop:deg_bound} that particular choices of the six points can give rise to irreducible plane quartics with three nodes in the exceptional locus. To make this example explicit consider $\calC$ an irreducible plane quartic with three nodes that we denote by $P_1,P_2$ and $P_6$. (Almost) any line through $P_2$ intersects the quartic in  two other points which we denote by $P_3$ and $R_1$, and similarly for the lines thorugh $P_6$, where the intersections we denote $P_5$ and $R_2$. Generically the intersection of the line $\ell_{23}$ through $P_2,P_3,R_1$ and the line $\ell_{56}$ through $P_5,P_6,R_2$ will not be a point of the quartic. But for a special choice of the points (still in general position) this can be made possible.
\end{example}

Examples \ref{ex:cubic}, \ref{ex:cubic2} and \ref{ex:quartic} show that, in certain cases, the exceptional locus can contain extra (singular) rational curves. It seems a difficult question to completely classify such exceptions: the problem is linked to the famous open problem to describe plane curves of fixed degree with given singular points (that goes back to results and questions of Severi).
%for \emph{every} configuration of the base points in general position the degree of $\pi(f(\PP^1))$ for a morphism $f: \PP^1 \to X$ with $f^*D$ supported in two points, is at most three, thus the previous examples form a complete list of exceptions.
%
%\begin{proposition}
%  Let $(X,D)$ be as in our setting and let $f: \PP^1 \to X$ be a morphism such that $f^*D$ is supported in at most two points. Then, either $f(\PP^1)$ is contracted by $\pi$ or $\deg \pi(f(\PP^1)) \leq 3$.
%\end{proposition}
\bibliography{references}{}

\begin{thebibliography}{RTW21b}

\bibitem[ABT20]{abt}
Kenneth Ascher, Lucas Braune, and Amos Turchet.
\newblock The {E}rd{\H o}s-{U}lam problem, {L}ang's conjecture and uniformity.
\newblock {\em Bull. Lond. Math. Soc.}, 52(6):1053--1063, 2020.

\bibitem[ADT20]{advt}
Kenneth Ascher, Kristin DeVleming, and Amos Turchet.
\newblock {Hyperbolicity and Uniformity of Varieties of Log General type}.
\newblock {\em International Mathematics Research Notices}, 08 2020.
\newblock online.

\bibitem[AT20]{AT_book}
Kenneth Ascher and Amos Turchet.
\newblock Hyperbolicity of varieties of log general type.
\newblock In Marc-Hubert Nicole, editor, {\em Arithmetic Geometry of
  Logarithmic Pairs and Hyperbolicity of Moduli Spaces: Hyperbolicity in
  Montr{\'e}al}, pages 197--247. Springer International Publishing, Cham, 2020.

\bibitem[Aut11]{Aut2}
Pascal Autissier.
\newblock Sur la non-densit\'{e} des points entiers.
\newblock {\em Duke Math. J.}, 158(1):13--27, 2011.

\bibitem[Bea96]{Beauville}
Arnaud Beauville.
\newblock {\em Complex algebraic surfaces}, volume~34 of {\em London
  Mathematical Society Student Texts}.
\newblock Cambridge University Press, Cambridge, second edition, 1996.
\newblock Translated from the 1978 French original by R. Barlow, with
  assistance from N. I. Shepherd-Barron and M. Reid.

\bibitem[BG06]{BG}
Enrico Bombieri and Walter Gubler.
\newblock {\em Heights in {D}iophantine geometry}, volume~4 of {\em New
  Mathematical Monographs}.
\newblock Cambridge University Press, Cambridge, 2006.

\bibitem[Cam04]{Ca04}
Fr\'{e}d\'{e}ric Campana.
\newblock Orbifolds, special varieties and classification theory.
\newblock {\em Ann. Inst. Fourier (Grenoble)}, 54(3):499--630, 2004.

\bibitem[CT21]{CaTur}
Laura Capuano and Amos Turchet.
\newblock Lang--{V}ojta conjecture over function fields for surfaces dominating
  {$\mathbb{G}^2_m$}.
\newblock {\em European Journal of Mathematics}, to appear, 2021.

\bibitem[CZ04]{CZAnnals}
Pietro Corvaja and Umberto Zannier.
\newblock On integral points on surfaces.
\newblock {\em Ann. of Math. (2)}, 160(2):705--726, 2004.

\bibitem[CZ08]{CZConic}
Pietro Corvaja and Umberto Zannier.
\newblock {Some cases of Vojta's Conjecture on integral points over function
  fields}.
\newblock {\em J. Algebraic Geometry}, 17:195--333, 2008.

\bibitem[CZ10]{CoZa10}
Pietro Corvaja and Umberto Zannier.
\newblock Integral points, divisibility between values of polynomials and
  entire curves on surfaces.
\newblock {\em Adv. Math.}, 225(2):1095--1118, 2010.

\bibitem[CZ13]{CZGm}
Pietro Corvaja and Umberto Zannier.
\newblock Algebraic hyperbolicity of ramified covers of {$\mathbb{G}^2_m$} (and
  integral points on affine subsets of {$\mathbb{P}_2$}).
\newblock {\em J. Differential Geom.}, 93(3):355--377, 2013.

\bibitem[DT15]{DT_invit}
Pranabesh Das and Amos Turchet.
\newblock Invitation to integral and rational points on curves and surfaces.
\newblock In {\em Rational points, rational curves, and entire holomorphic
  curves on projective varieties}, volume 654 of {\em Contemp. Math.}, pages
  53--73. Amer. Math. Soc., Providence, RI, 2015.

\bibitem[Fal91]{Falt_Ab}
Gerd Faltings.
\newblock Diophantine approximation on abelian varieties.
\newblock {\em Ann. of Math. (2)}, 133(3):549--576, 1991.

\bibitem[GSW21]{GSW}
Ji~Guo, Chia-Liang Sun, and Julie Tzu-Yueh Wang.
\newblock Some cases of {V}ojta's conjectures related to algebraic tori over
  function fields.
\newblock {\em arXiv preprint arXiv:2106.15881}, 2021.

\bibitem[Har77]{Hartshorne}
Robin Hartshorne.
\newblock {\em Algebraic geometry}.
\newblock Springer-Verlag, New York-Heidelberg, 1977.
\newblock Graduate Texts in Mathematics, No. 52.

\bibitem[HS00]{SH}
Marc Hindry and Joseph~H. Silverman.
\newblock {\em Diophantine geometry}, volume 201 of {\em Graduate Texts in
  Mathematics}.
\newblock Springer-Verlag, New York, 2000.

\bibitem[Jav20]{J_book}
Ariyan Javanpeykar.
\newblock The lang--vojta conjectures on projective pseudo-hyperbolic
  varieties.
\newblock In Marc-Hubert Nicole, editor, {\em Arithmetic Geometry of
  Logarithmic Pairs and Hyperbolicity of Moduli Spaces: Hyperbolicity in
  Montr{\'e}al}, pages 135--196. Springer International Publishing, Cham, 2020.

\bibitem[JK20]{JK}
Ariyan Javanpeykar and Ljudmila Kamenova.
\newblock Demailly's notion of algebraic hyperbolicity: geometricity,
  boundedness, moduli of maps.
\newblock {\em Math. Z.}, 296(3-4):1645--1672, 2020.

\bibitem[Lan86]{Lan86}
Serge Lang.
\newblock Hyperbolic and {D}iophantine analysis.
\newblock {\em Bull. Amer. Math. Soc. (N.S.)}, 14(2):159--205, 1986.

\bibitem[Lan91]{Lang91}
Serge Lang.
\newblock {\em Number theory. {III} - Diophantine geometry}, volume~60 of {\em
  Encyclopaedia of Mathematical Sciences}.
\newblock Springer-Verlag, Berlin, 1991.

\bibitem[Lev09]{Levin}
Aaron Levin.
\newblock Generalizations of {S}iegel's and {P}icard's theorems.
\newblock {\em Ann. of Math. (2)}, 170(2):609--655, 2009.

\bibitem[McQ98]{McQ98}
Michael McQuillan.
\newblock Diophantine approximations and foliations.
\newblock {\em Inst. Hautes \'Etudes Sci. Publ. Math.}, No. 87:121--174, 1998.

\bibitem[RTW21a]{RTW2}
Erwan Rousseau, Amos Turchet, and Julie Tzu-Yueh Wang.
\newblock Divisibility of polynomials and degeneracy of integral points.
\newblock {\em arXiv preprint arXiv:2106.11337}, 2021.

\bibitem[RTW21b]{RTW}
Erwan Rousseau, Amos Turchet, and Julie Tzu-Yueh Wang.
\newblock Nonspecial varieties and generalised {L}ang-{V}ojta conjectures.
\newblock {\em Forum Math. Sigma}, 9:e11, 29, 2021.

\bibitem[RV20]{RV19}
Min Ru and Paul Vojta.
\newblock A birational {N}evanlinna constant and its consequences.
\newblock {\em Amer. J. Math.}, 142(3):957--991, 2020.

\bibitem[Tur17]{Tur}
Amos Turchet.
\newblock Fibered threefolds and {L}ang-{V}ojta's conjecture over function
  fields.
\newblock {\em Trans. Amer. Math. Soc.}, 369(12):8537--8558, 2017.

\bibitem[Voj87]{Vojta87}
Paul Vojta.
\newblock {\em {Diophantine Approximations and Value Distribution Theory}},
  volume 1239 of {\em Lecture Notes in Mathematics}.
\newblock Springer Berlin Heidelberg, 1987.

\bibitem[Voj96]{VojtaSA}
Paul Vojta.
\newblock Integral points on subvarieties of semiabelian varieties. {I}.
\newblock {\em Invent. Math.}, 126(1):133--181, 1996.

\bibitem[Voj99]{VojtaSA2}
Paul Vojta.
\newblock Integral points on subvarieties of semiabelian varieties. {II}.
\newblock {\em Amer. J. Math.}, 121(2):283--313, 1999.

\bibitem[Voj11]{Vojta11}
Paul Vojta.
\newblock Diophantine approximation and {N}evanlinna theory.
\newblock In {\em Arithmetic geometry}, volume 2009 of {\em Lecture Notes in
  Math.}, pages 111--224. Springer, Berlin, 2011.

\bibitem[Wal78]{Walker}
Robert~J. Walker.
\newblock {\em Algebraic curves}.
\newblock Springer-Verlag, New York-Heidelberg, 1978.
\newblock Reprint of the 1950 edition.

\end{thebibliography}
\bibliographystyle{alpha}

\end{document}